\documentclass{amsart}
\usepackage{amsmath,amsfonts,epic,eepic}
\usepackage[all,cmtip]{xy}

\newtheorem{theorem}{Theorem}
\newtheorem{lemma}{Lemma}
\newtheorem{definition}{Definition}

\newcommand{\im}{\operatorname{im}}
\newcommand{\ord}{\operatorname{ord}}

\newcommand{\Z}{\mathbb{Z}}

\newtheorem{corollary}{Corollary}

\title{Smoothness theorem for differential BV algebras}

\author{John Terilla}
\address{Mathematics Department, CUNY Graduate Center and 
 Queens College, City University of
  New York, 65-30 Kissena Blvd, Flushing NY 11367, USA}
  \email{jterilla@qc.cuny.edu} 

\thanks{Supported in part by a City
    University of New York Collaborative Research Grant.}


\begin{document}

\begin{abstract}
  Given a differential Batalin-Vilkovisky algebra $(V,Q,\Delta,
  \cdot)$, the associated odd differential graded Lie algebra
  $(V,Q+\Delta, [\,,\,])$ is always smooth formal.  The 
  quantum dgLa   $L_\hbar:=\left(V[[\hbar]], Q+\hbar \Delta, [\,,\,]\right)$ is not always
  smooth formal, but when it is---for example when a $Q$-$\Delta$ version
  of the $\partial$-$\overline{\partial}$ Lemma holds---there is a
  weak-Frobenius manifold structure on the homology of $L$ that is
  important in applications and relevant to quantum correlation
  functions.  In this paper, we prove that $L_\hbar$ is smooth formal
  if and only if the spectral sequence associated to the filtration
  $F^p:=\hbar^p V[[\hbar]]$ on the complex $(V[[\hbar]],Q+\hbar
  \Delta)$ degenerates at $E_1$.  A priori, this degeneration means
  that a collection of first order obstructions vanish and we prove
  that it follows that all obstructions vanish.  For those
  differential BV algebras that arise from the Hochschild complex of a
  Calabi-Yau category, the degeneration of this spectral sequence is
  an expression of the noncommutative Hodge to deRham degeneration,
  conjectured by Kontsevich and Soibelman and proved to hold in
  certain cases by   Kaledin.  
  The results in this paper imply that the noncommutative Hodge to de Rham
  degeneration conjecture is equivalent to the existence of a versal
  solution to the quantum master equation.  
  At the end of the paper, some physical considerations are mentioned.
\end{abstract}

\maketitle
\section{Introduction}
An odd differential ($\Z$ or $\Z/2$) graded Lie algebra
$(V,Q,[\,,\,])$ is \emph{smooth formal} if and only if there exists a
versal solution to the Maurer-Cartan equation; that is, a degree zero
solution
\[\Gamma=\sum \gamma_it^i+\gamma_{ij}t^i t^j+\gamma_{ijk}t^i t^j
t^k+\cdots \] to the equation \[Q\Gamma+\frac{1}{2}[\Gamma,\Gamma]=0\]
where the power series coordinates $\{t^i\}$ are a homogeneous basis
for the dual of the graded vector space $H:=\ker(Q)/\im(Q)$ and
$\{\gamma_i\}$ are representatives for homology classes dual to the
$\{t^i\}$.  The terminology ``smooth formal'' arises from the fact
that a differential graded Lie algebra has a versal solution to the
Maurer-Cartan equation if and only if it is quasi-isomorphic to its
homology (i.e. formal) with zero bracket (i.e. smooth).  Equivalently,
a differential graded Lie algebra is smooth formal if and only if it
represents a smooth point in its Maurer-Cartan moduli space, or, if
its associated deformation functor is a smooth functor.

\begin{definition}A Batalin-Vilkovisky (BV) algebra is a triple
  $(V,\Delta,\cdot)$ where $(V,\cdot)$ is a $\Z$ graded, commutative,
  associative, unital algebra and $\Delta:V^k \to V^{k-1}$ satsifies
  $\Delta^2=0$, and for each $v\in V^k$, the operator $d_v$ defined by
  $d_v(w)=\Delta (v\cdot w)-\Delta(v)\cdot w-(-1)^{k}v\cdot \Delta(w)$
  is a derivation of degree $k-1$.
\end{definition}

We suppress the symbol for the associative product and denote the BV
algebra $(V,\Delta, \cdot)$ simply by $(V,\Delta)$ and denote $v\cdot
w$ by $vw$.

Given a BV algebra $(V,\Delta)$, the operation $[\,,]:V^k\otimes V^l
\to V^{k+l-1}$ defined by $[v,w]:=d_v(w)$ satisfies:
\begin{align*}
  [v,w]&=-(-1)^{(i+1)(j+1)}[w,v]\\
  [v,[w,u]]&=[[v,w],u]-(-1)^{(i+1)(j+1)}[w,[v,u]] \intertext{and}
  \Delta[v,w]&=[\Delta(v),w]+(-1)^{i+1}[v,\Delta(w)]
\end{align*}
for $v\in V^i$, $w\in V^j$, and $u\in V$.  Thus, the triple
$(V,\Delta, [\,,\,])$ is an odd differential graded Lie algebra.
However, up to quasi-isomorphism, it is not a very interesting odd
dgLa as the following theorem establishes.

\begin{theorem}\label{class_smooth_thm}
  Let $(V,\Delta)$ be a BV algebra and let $L=(V,\Delta, [\,,\,])$ be
  the associated odd dgLa.  Then $L$ is smooth formal.
\end{theorem}

We know several ways to prove this theorem.  It is straightforward to
produce a versal $\Gamma$ term by term.  Begin with $\Gamma_1:=\sum
\gamma_i t^i$ where $\{[\gamma_i]\}$ is a homogeneous basis of $H$
dual to $\{t^i\}$.  Then, the quadratic term $\Gamma_2=\sum
\gamma_{ij}t^it^j$ must satisfy \[\Delta \Gamma_2 +
\frac{1}{2}[\Gamma_1, \Gamma_1]=0\Leftrightarrow \Delta
\gamma_{ij}=-\frac{1}{2}[\gamma_i,\gamma_j].\] One solution for $\gamma_{ij}$ is
$\gamma_{ij}=-\frac{1}{2}\gamma_i \gamma_j$ since
\[[\gamma_i,\gamma_j]=\Delta(\gamma_i
\gamma_j)-\Delta(\gamma_i)\gamma_j-(-1)^{|\gamma_i|} \gamma_i
\Delta(\gamma_j)=\Delta(\gamma_i\gamma_j).\] One solves $\Delta
\Gamma_3+\frac{1}{2}[\Gamma_1,\Gamma_2]+\frac{1}{2}[\Gamma_2,\Gamma_1]=0
\Leftrightarrow \Delta \Gamma_3 = -[\Gamma_1,\Gamma_2]$ for the cubic
term $\Gamma_3=\sum \gamma_{ijk}t^it^jt^k$ and finds a solution for
$\gamma_{ijk}$ in terms of products of $\gamma_{i}$, $\gamma_{j}$, and
$\gamma_{k}$, and so on.  One finds that the cohomological
obstruction to finding $\Gamma_n$ vanishes and one can solve for 
$\gamma_{i_1 \ldots i_n}$ by induction.

There is a conceptually better proof based on the following Lemma:
\begin{lemma}\label{lemma1}
  If $(V,\Delta)$ is a BV algebra and $\gamma\in V^0$ then
  \[\Delta(\gamma)+\frac{1}{2}[\gamma, \gamma]=0 \Leftrightarrow
  \Delta \left(e^{\gamma}\right)=0.\]
\end{lemma}
\begin{proof}
  For $\gamma \in V^0$, one has
  $\Delta(\gamma^2)=2\gamma\Delta(\gamma) + [\gamma, \gamma]$ and by
  induction
  \[\Delta(\gamma^n)=n\gamma^{n-1}\Delta(\gamma)+\binom{n}{2}\gamma^{n-2} [\gamma, \gamma]\] from which it follows
  \begin{align*}
    \Delta(e^\gamma)&=\Delta
    \left ( \sum_{n=0}^\infty \frac{1}{n!}\gamma^n \right)\\
    &=\Delta(\gamma)\sum_{n=1}^\infty
    \frac{1}{(n-1)!}\gamma^{n-1}+\frac{1}{2}[\gamma,
    \gamma]\sum_{n=2}^\infty \frac{1}{(n-2)!}\gamma^{n-2}\\
    &=\left(\Delta(\gamma)+\frac{1}{2}[\gamma, \gamma]\right)e^\gamma.
  \end{align*}
  Since $e^\gamma$ is a unit, $\Delta(e^\gamma)=0 \Leftrightarrow
  \Delta(\gamma)+\frac{1}{2}[\gamma, \gamma]=0.$
\end{proof}
Thus, to prove Theorem \ref{class_smooth_thm}, let $\{\gamma_i\}$ be
representatives for a basis for $H$, set \[\Gamma=\log\left(1+\sum_i
  \gamma_i t^i\right)=\sum_i\gamma_i t^i -
\frac{1}{2}\sum_{ij}\gamma_i\gamma_j t^i t^j + \frac{1}{3}
\sum_{ijk}\gamma_i\gamma_j \gamma_k t^i t^j -\cdots\] from which it
follows that \[\Delta\left(e^\Gamma \right)=\Delta\left(1+\sum
  \gamma_i t^i\right)=0.\]

So, if an odd dgLa originates from BV data, it is smooth formal---the
Lie bracket vanishes in homology and so do all the higher Lie-Massey
products, which are the obstructions to solving the Maurer-Cartan
equation.

\section{Differential BV algebras and the main theorem}

\begin{definition}A differential BV (dBV) algebra is a fourtuple
  $(V,Q, \Delta,\cdot)$ where $(V,Q, \cdot)$ is a differential graded,
  commutative, associative, unital algebra, $(V,\Delta,
  \cdot)$ is a BV algebra, and $Q\Delta+\Delta Q=0$.
\end{definition}
Again, we suppress the notation for the commutative associative
product and write the data of a differential BV algebra as $(V, Q,
\Delta)$.  In this paper, we adopt the
convention that $Q$ has degree $+1$ and $\Delta$ has degree $-1$,
since these signs match the naturally occuring 
example of the Hochschild cochains
of a Calabi-Yau category.
However all results hold when the degrees of $Q$ and $\Delta$ 
are either $\pm 1$, either the same or different, or if $V$ is only
$\Z/2$ graded.

Given a dBV algebra $(V,Q,\Delta)$, one considers
the bracket defined 
as usual using the BV operator $\Delta$:
\begin{align*} [v,w]:
=\Delta (v\cdot w)-\Delta(v)\cdot w-(-1)^{|v|}v\cdot
  \Delta(w).
\end{align*}
There are several related odd differential
graded Lie algebras to consider.  Disregarding $(V,\Delta, [\,,\,])$
and $(V,Q+\Delta, [\,,\,])$
which are always smooth formal by Theorem
\ref{class_smooth_thm}, 
let
\begin{align}
  L&:=(V, Q, [\,,\,])\\
  L_\hbar&:=(V[[\hbar]], Q+\hbar\Delta, [\,,\,])
\end{align}
We call the first one the \emph{classical dgLa} and the second one the
\emph{quantum dgLa} associated to the dBV algebra $(V,Q,\Delta)$.  The
quantum dgLa $L_\hbar$ is a Lie algebra over the ring $k[[\hbar]]$
where $k$ is the ground field of the vector space $V$, and $\hbar$ is
a formal variable, which, using our degree conventions for $Q$ and
$\Delta$, is assigned the 
degree $+2$.  For
convenience, let
\begin{equation}K:=Q+\hbar \Delta.
\end{equation}
We say the dgLa $L_\hbar$ is \emph{smooth formal as a quantum dgLa}
if and only if there exists a versal solution to the quantum master
equation; that is, if and only if there exists
\begin{equation}\label{Quantum_Gamma}
\Gamma=\sum \gamma_it^i+\gamma_{ij}t^i t^j+\gamma_{ijk}t^i
t^j t^k+\cdots
\end{equation} satisfying
\begin{equation}\label{QME}
  K\Gamma+\frac{1}{2}[\Gamma,\Gamma]=0
\end{equation} 
where the power series coordinates $\{t^i\}$ comprise
a homogeneous basis for the dual of the graded vector space
$H:=\ker(Q)/\im(Q)$,
the terms $\gamma_{i_1 \cdots i_n}$ are expanded
in powers of $\hbar$ as
\[\gamma_{i_1\cdots
  i_n}=\gamma_{i_1\cdots i_n}^{(0)}+\gamma_{i_1\cdots
  i_n}^{(1)}\hbar+\gamma_{i_1\cdots i_n}^{(2)}\hbar^2+\cdots \in
V[[\hbar]],\]
and $\left\{\gamma_i^{(0)}\right\}$ are
representatives of homology classes dual to the $\{t^i\}$.  Note that
the smooth formality of the quantum Lie algebra $L_\hbar$ implies the
smooth formality of the classical Lie algebra $L$ since a versal 
solution to the quantum master equation
$K\Gamma+\frac{1}{2}[\Gamma,\Gamma]=0$ maps to a solution to the
Maurer-Cartan equation $Q\Gamma+\frac{1}{2}[\Gamma,\Gamma]=0$ 
under the projection
$\alpha:V[[\hbar]]\to V$ defined by setting $\hbar=0$.

Unlike $(V,\Delta,[\,,\,])$ and 
$(V,Q+\Delta,[\,,\,])$, the smooth formality of the quantum dgLa
$L_\hbar=(V[[\hbar]],Q+\hbar\Delta, [\,,\,])$ and the classical dgLa
$L=(V,Q,[\,,\,])$ 
are not automatic, and so it is often an interesting situation  
when they are smooth formal.  For $v\in
V^k[[\hbar]]$ and $w\in V[[\hbar]]$, the relationship between the
associative product, the Lie bracket, and differential $K$ is
expressed
\[K(vw)+K(v)w+(-1)^kvK(w)=\hbar[v,w].\]  
So, if $\gamma_i$ and $\gamma_j$ are both
$K$-closed, then $\hbar[\gamma_i,\gamma_j]=K(\gamma_i\gamma_j)$ is
$K$-exact, but $[\gamma_i,\gamma_j]$ may not be $K$-exact in $V[[\hbar]]$---the simple trick that proves
that $(V,\Delta,[\,,\,])$ and 
$(V,Q+\Delta,[\,,\,])$ are smooth formal does not apply in
 $L_\hbar$.  But there are even more fundamental
obstructions to the smooth formality of $L_\hbar$.  The
initial condition of beginning with 
a basis $\{\gamma_i^{(0)}\}$ of representatives for the classical homology
$H=H(V,Q)$ provides an infinitesimal solution
$\Gamma_1^{(0)}:=\sum_i\gamma_i^{(0)}t^i$ to the master equation modulo
$t^it^j$ and modulo $\hbar$.  To build the solution $\Gamma$ term by
term, picture
the terms arranged as in Figure 1
in the first quadrant with the
$x$-axis labelled by the powers of $\hbar$ and the $y$-axis labelled
by the powers of $t$.  The infinitesimal solution $\Gamma_1^{(0)}$
resides at the $(1,1)$ position.  In order to extend $\Gamma_1^{(0)}$ 
 to the right---the positions circled in the picture---by adding terms with higher powers of $\hbar$ so that
the result
\[\Gamma_1^{(0)}+\Gamma_1^{(1)}+\cdots+\Gamma_1^{(r)}=
\sum_i \left(\gamma_i^{(0)}+
\gamma_i^{(1)}\hbar+\cdots+\gamma_i^{(r)}\hbar^r\right)t^i,\]
is a solution to the quantum master equation modulo $t^it^j$ and
modulo $\hbar^{r+1}$
it must be that $K\left(\Gamma_1^{(0)}+\Gamma_1^{(1)}+\cdots +\Gamma_1^{(r)}\right)=0$
modulo $\hbar^{r+1}$
which means that
$-\Delta\left(\gamma_i^{(j)}\right)=Q\left(\gamma_i^{(j+1)}\right)$
for $j=0,\ldots, r-1$.  

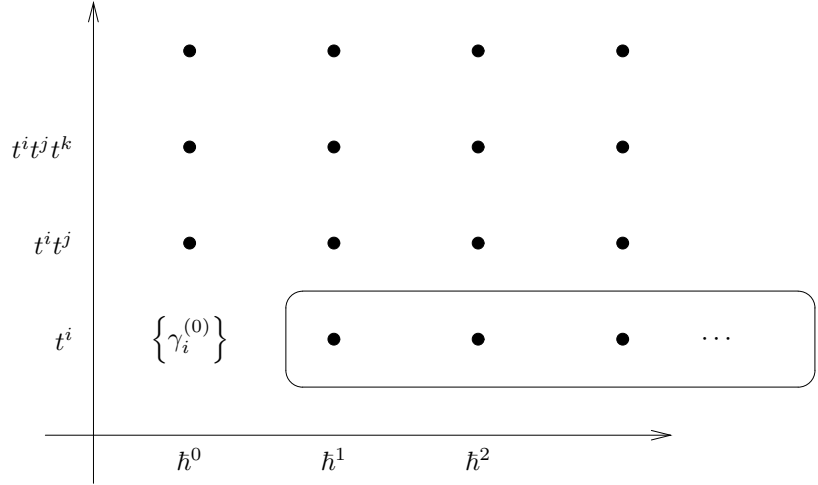
\begin{figure}[htb]\label{array_of_obstructions}
\caption{The array of obstructions}
\begin{center}
\setlength{\unitlength}{0.00083333in}
{\renewcommand{\dashlinestretch}{30}
\begin{picture}(3987,3039)(0,-10)
\put(1680,717){\arc{210}{1.5708}{3.1416}}
\put(1680,1107){\arc{210}{3.1416}{4.7124}}
\put(4770,1107){\arc{210}{4.7124}{6.2832}}
\put(4770,717){\arc{210}{0}{1.5708}}
\path(1575,717)(1575,1107)
\path(1680,1212)(4770,1212)
\path(4875,1107)(4875,717)
\path(4770,612)(1680,612)
\put(4275,912){\makebox(0,0)[c]{$\cdots$}}
\put(975,912){\makebox(0,0)[c]{$\left\{\gamma_i^{(0)}\right\}$}}
\path(375,12)(375,3012)
\path(405.000,2892.000)(375.000,3012.000)(345.000,2892.000)
\put(250,912){\makebox(0,0)[r]{$t^i$}}
\put(250,1512){\makebox(0,0)[r]{$t^it^j$}}
\put(250,2112){\makebox(0,0)[r]{$t^it^jt^k$}}
\put(975,162){\makebox(0,0)[c]{$\hbar^0$}}
\put(1875,162){\makebox(0,0)[c]{$\hbar^1$}}
\put(2775,162){\makebox(0,0)[c]{$\hbar^2$}}
\path(75,312)(3975,312)
\path(3855.000,282.000)(3975.000,312.000)(3855.000,342.000)
\put(975,1512){\blacken\ellipse{75}{75}}
\put(975,2112){\blacken\ellipse{75}{75}}
\put(1875,912){\blacken\ellipse{75}{75}}
\put(1875,1512){\blacken\ellipse{75}{75}}
\put(1875,2112){\blacken\ellipse{75}{75}}
\put(1875,2712){\blacken\ellipse{75}{75}}
\put(975,2712){\blacken\ellipse{75}{75}}
\put(2775,912){\blacken\ellipse{75}{75}}
\put(2775,1512){\blacken\ellipse{75}{75}}
\put(2775,2112){\blacken\ellipse{75}{75}}
\put(2775,2712){\blacken\ellipse{75}{75}}
\put(3675,2712){\blacken\ellipse{75}{75}}
\put(3675,2112){\blacken\ellipse{75}{75}}
\put(3675,1512){\blacken\ellipse{75}{75}}
\put(3675,912){\blacken\ellipse{75}{75}}
\end{picture}
}
\end{center}
\end{figure}

\subsection{The $\hbar$ filtration and associated spectral sequence}
Forget for a moment the associative and the Lie products, and consider the
complex $(V[[\hbar]],K)$.  It is decreasingly filtered by powers of $\hbar$:
\[V[[\hbar]]=F^0 \supset F^1\supset F^2 \supset \cdots \text{ where
}F^p:=\hbar^p V[[\hbar]]\] Let $(E_r^p,d_r)$ be the associated
spectral sequence.  Note that $k[[\hbar]]$ acts on the spectral
sequence in a simple way \[\hbar E^p_r=E^{p+1}_r\] so the entire
spectral sequence is determined by $E^0_r$.  The differential
$K=Q+\hbar\Delta$ acts on $E_0^0=F^0/F^1$ by $Q$ and so $E_1^0$ is
precisely the classical cohomology $H=H(V,Q)$.  A classical class in
$E_1^0$ survives to $E^0_\infty$ provided it can be extended to give a
quantum class; i.e., provided there exists a representative
$\gamma^{(0)}$ of the classical cohomology class and elements
$\gamma^{(r)}\in V$ so that $K(\gamma^{(0)}+\hbar \gamma^{(1)}+\hbar^2
\gamma^{(2)}+\cdots)=0.$ The spectral sequence
degenerates at the $E_1$ term if and only if there is a basis
$\left\{\gamma_i^{(0)}\right\}$ of $H(V,Q)$ and elements
$\gamma_i^{(j)}\in V$ so that
\[K\left(\gamma_i^{(0)}+\hbar \gamma_i^{(1)}+\hbar^2
  \gamma_i^{(2)}+\cdots\right)=0.\] The degeneration of this spectral
sequence was conjectured by Kontsevich and Soibelman \cite{KS} to hold
for the Hochschild complex of a Calabi-Yau category and was proved to
hold in cases by Kaledin \cite{Ka}.  We shall refer to the
degeneration of $(E_r^p,dr)$ at the $E_1$ term simply as \emph{the
  degeneration of the spectral sequence}.

Now, returning to the master equation and the array of obstructions in
Figure 1, we see that the degeneration of the spectral sequence is
equivalent to the vanishing of the first row of obstructions, the ones
circled in Figure 1.  So, the degeneration of the spectral sequence
is necessary for the smooth formality of $L_\hbar$.  The next theorem,
and the main result of the paper, asserts that it is also sufficient.

\begin{theorem}\label{main_thm}
  $L_\hbar$ is smooth formal as a quantum dgLa if and only if the
  spectral sequence associated to the filtration
  $F^p=\hbar^pV[[\hbar]]$ of the complex $(V[[\hbar]],K)$ degenerates
  at the $E_1$ term.
\end{theorem}
\begin{proof}
  The map $V[[\hbar]]\to V$ defined by setting $\hbar=0$ is a chain
  map, denote it by $\alpha$.  The degeneration of the spectral
  sequence means that there exists a chain complex splitting,
  denote it by $\beta$:
  \[\xymatrix{(V[[\hbar]],K) \ar[r]_-\alpha & (V,Q) \ar@/_1pc/[l]_-\beta
    \ar[r] & 0}\]
  We use the maps $\alpha$ and $\beta$ to construct a solution
  \[\Gamma=\sum \gamma_it^i+\gamma_{ij}t^i t^j+\gamma_{ijk}t^i
  t^j t^k+\cdots\] to the quantum master equation.  We work with all
  powers of $\hbar$ at once, and the terms
  \[\Gamma_k:= \sum \gamma_{i_1i_2\cdots i_k}t^{i_1}t^{i_2}\cdots
  t^{i_k}\] will be defined inductively so that $\Gamma_1+\cdots
  +\Gamma_k$ is a solution to the equation
  \[K\Gamma+\frac{1}{2}[\Gamma,\Gamma]=0 \mod t^{i_1}\cdots
  t^{i_{k+1}}.\]
  The first step is to let $\left\{\gamma_i^{(0)}\right\}$ be a basis
  of representatives for $H=H(V,Q)$ and define \[\Gamma_1=\sum_i
  \beta\left(\gamma_i^{(0)}\right) t^i.\]  Note that $K\beta
  \gamma_i^{(0)}=\beta Q \gamma_i^{(0)}=0$ so
  $\Gamma_1$ is a solution to the master equation modulo $t^it^j$.

  We perform the second step explicitly.  Being able to extend
  $\Gamma_1$ to $\Gamma_1+\Gamma_2$ so that the result is a solution
  to the master equation
  modulo $t^it^jt^k$ is equivalent to
  $-\frac{1}{2}[\Gamma_1,\Gamma_1]$ being $K$-exact. Note that
  $\hbar [\Gamma_1,\Gamma_1]$ is $K$-exact since
  $K(\Gamma_1^2)=K(\Gamma_1\Gamma_1)=\hbar[\Gamma_1,\Gamma_1]$.  We
  make a correction by subtracting a $K$-closed term from $\Gamma_1^2$ so
  that the difference is divisible by $\hbar$.  Note
  \begin{gather}
    K\beta \alpha \Gamma_1^2=\beta Q \alpha \Gamma_1^2=\beta \alpha
    K\Gamma_1^2=\beta \alpha \hbar [\Gamma_1,\Gamma_1]=0\label{closed}
    \intertext{and} \alpha(\Gamma_1^2-\beta \alpha \Gamma_1^2)=\alpha
    \Gamma_1^2-\alpha \beta \alpha \Gamma_1^2 =\alpha
    \Gamma_1^2-\alpha \Gamma_1^2=0.
    \label{divisible}
  \end{gather}
  Equation \eqref{closed} shows that $\beta\alpha\Gamma_1^2$ is
  $K$-closed and Equation \eqref{divisible} shows that the difference
  $\Gamma_1^2-\beta \alpha \Gamma_1^2$ is divisible by $\hbar$.
  Therefore, $K\left(\frac{1}{\hbar}\left(\Gamma_1^2-\beta \alpha
      \Gamma_1^2\right)\right)=[\Gamma_1,\Gamma_1]$ and setting
  \[\Gamma_2:=-\frac{1}{2\hbar}\left(\Gamma_1^2-\beta \alpha
    \Gamma_1^2\right)
  \]
  gives
  \[K(\Gamma_1+\Gamma_2)+\frac{1}{2}[\Gamma_1+\Gamma_2,\Gamma_1+\Gamma_2]=0
  \mod t^it^jt^k.
  \]
  
For the general inductive step, observe that the reasoning used in Equations \eqref{closed} and \eqref{divisible}  
  applies generally and yields \eqref{closed_general} and
  \eqref{divisible_general}:
\begin{gather}
Ky\in \hbar V[[\hbar]]\Rightarrow K\beta \alpha y =0 \label{closed_general}
\intertext{and}
y\in V[[\hbar]]\Rightarrow (1-\beta \alpha)y\in \hbar
    V[[\hbar]].\label{divisible_general}
\end{gather}


  Now suppose that $\Gamma_1, \ldots, \Gamma_{n-1}$ have been defined
  so that $\Gamma_1+\cdots + \Gamma_{n-1}$ satisfies
  \[K(\Gamma_1+\cdots + \Gamma_{n-1})+\frac{1}{2}[\Gamma_1+\cdots +
  \Gamma_{n-1},\Gamma_1+\cdots + \Gamma_{n-1}]=0 \mod t^{i_1}\cdots
  t^{i_n}.\] Let 
  \begin{align*}
    x&=\ord_n\left(K(\Gamma_1+\cdots +
  \Gamma_{n-1})+\frac{1}{2}[\Gamma_1+\cdots +
  \Gamma_{n-1},\Gamma_1+\cdots + \Gamma_{n-1}]\right)\\
&= \ord_n\left(\frac{1}{2}[\Gamma_1+\cdots +
  \Gamma_{n-1},\Gamma_1+\cdots + \Gamma_{n-1}]\right)
\end{align*} where $\ord_n$ means the
term of order $n$ in the variables $\{t^i\}$.   We will show that $x$
is $K$-exact.   Then, for $\Gamma_n\in V[[\hbar]]$ 
  satisfying $K(\Gamma_n)=-x$, the sum
  $\Gamma_1+\cdots +\Gamma_{n-1}+\Gamma_n$ satisfies the master
  equation modulo $t^{i_1}\cdots t^{i_{n+1}}$ as needed.

  Following the steps in the calculation in Lemma
  \ref{lemma1} and keeping track of the $\frac{1}{\hbar}$ shows that
  for any any degree zero $\gamma\in V[[\hbar]][t^i]$,
     \[e^{-\frac{\gamma}{\hbar}}\hbar K\left(
       e^{\frac{\gamma}{\hbar}}\right)=K(\gamma)+\frac{1}{2}[\gamma,
     \gamma]. \]
  In particular, we can rewrite $x=\ord_n\left(e^{-\frac{\Gamma_1+\cdots +
      \Gamma_{n-1}}{\hbar}}\hbar K\left( e^{\frac{\Gamma_1+\cdots +
        \Gamma_{n-1}}{\hbar}}\right)\right).$  Multiplying by 
\[e^\frac{\Gamma_1+\cdots +
      \Gamma_{n-1}}{\hbar}=1+\left(\Gamma_1+\cdots
      +\Gamma_{n-1}\right)+\cdots\] only affects terms of order $n+1$
    and higher, hence 
\begin{equation}\label{x_rewritten}
x=\ord_n\left(\hbar K\left(e^{\frac{\Gamma_1+\cdots +
          \Gamma_{n-1}}{\hbar}}\right)\right).
\end{equation}
This doesn't show that $x$ is $K$-exact since $\left(e^{\frac{\Gamma_1+\cdots +
          \Gamma_{n-1}}{\hbar}}\right)$ has negative powers of $\hbar$.
  But multiplying Equation \eqref{x_rewritten} by $\hbar^{n-1}$ shows that $\hbar^{n-1}x$ is
  $K$-exact since 
\[\hbar^{n-1}x=Ky \text{ where }y=\ord_n\left(\hbar^n e^{\frac{\Gamma_1+\cdots +
          \Gamma_{n-1}}{\hbar}}\right)\] and $y$ has no negative
    powers of $\hbar$.
  Then the implications in \eqref{closed_general} and \eqref{divisible_general} 
 show that $K(\beta \alpha y)=0$ and $(1-\beta\alpha)y$ 
 is divisible by $\hbar$.
  Therefore,
  \[y_{n-1}:=\frac{1}{\hbar}\left(1 -\beta \alpha \right)y \]
  satisfies $K(y_{n-1})=\hbar^{n-2}x.$ Iterating, we find for $k=1,
  \ldots, n-1$
  \[y_{n-k}:=\left(\frac{1 -\beta \alpha}{\hbar}\right)^{k} y
  \]
  satisfies $K(y_{n-k})=\hbar^{n-k-1}x.$ The outcome is a term $y_1$,
  with no negative powers of $\hbar$, satisfying
  $K(y_1)=x$ as needed.

\end{proof}
\section{Examples}

\subsection{Differential BV algebras satisfying the $Q$-$\Delta$ Lemma}
In \cite{BK,M} it is proved that if $(V,Q,\Delta)$ is a differential
BV algebra that satisfies the $Q$-$\Delta$ version of the
$\partial$-$\bar{\partial}$ lemma, then the classical dgLa
$(V,Q,[\,,\,])$ is smooth formal.  
Recall a complex with two differentials $Q$ and $\Delta$ is said to
satisfy the $Q$-$\Delta$ lemma provided
\begin{equation}
\im(Q\Delta)=\im(\Delta Q)=\im(Q)\cap \ker(\Delta)=\im(\Delta)\cap
\ker(\Delta).
\end{equation}
One consequence of the $Q-\Delta$ lemma is that every $Q$ cohomology
class has a representative that is $\Delta$ closed, hence $K$ closed
in $V[[\hbar]]$.  No ``$\hbar$ corrections'' are necessary.
This gives a splitting $\beta:V\to V[[\hbar]]$ (of a type stronger
than necessary)
and consequently, the spectral sequence degenerates.  
Therefore, we have the following consequence of Theorem \ref{main_thm}.
\begin{corollary}
If $(V,Q,\Delta)$ is a differential $BV$ algebra satisfying the
$Q$-$\Delta$ lemma, then $L_\hbar$ is smooth formal over $k[[\hbar]]$,
hence $L$ is smooth formal.
\end{corollary}


\subsection{A polynomial algebra example that doesn't satisfy the
  $Q$-$\Delta$ Lemma}
Singularity theory and Landau-Ginzberg models 
provides many examples of dBV algebras that do not
satisfy a $Q$-$\Delta$ lemma, but do give rise to quantum dgLa's that
are smooth formal.  Here is a simple example.

Let $U=U^{-1}\oplus U^{0}$ be a two dimensional vector space with
basis $\{\eta,x\}$.  We define a BV algebra $(V,\Delta)$ as follows:
$V=SU$ is the free commutative associative algebra on $U$ and
\[\Delta=\frac{\partial ^2}{\partial x\partial \eta}.\] The bracket
defined by $\Delta$ is seen to be the one on $V$ given on the
generators by
 \[[x,x]=[\eta,\eta]=0 \text{ and }[x,\eta]=-[\eta,x]=1\]
 and extended to $V$ as a derivation of the associative product in
 each coordinate.  For analogy, $U$ is a graded symplectic vector
 space and $V=SU$ is like the Poisson algebra of (polynomial)
 functions on $U$.  Let $S^{(0)}=x^3$ and note that
 it satisfies $[S^{(0)},S^{(0)}]=0$ and $\Delta(S^{(0)})=0$ hence
 defines a solution to the quantum master equation.  We have the
 differential $Q$ given by $Q(a)=[S,a]$.  We have defined a 
differential BV algebra $(V,Q, \Delta)$ and have the two dgLa's
 $L=(V,Q,[\,,\,])$ and $L_\hbar =(V[[\hbar]],K,[\,,\,]).$ Explicitly,
 \begin{gather*}
   L=L^{-1}\oplus L^0 \text{ where } L^{-1}=k[x]\eta \text{ and
   }L^0=k[x],\\
   Q(f+g\eta) = 3x^2 g \text{ for }f,g \in k[x],\\
   \Delta(f+g\eta) = \frac{\partial g}{\partial x} \text{ for }f,g \in k[x],\\
   K(f+g\eta)=3x^2g+\hbar\frac{\partial g}{\partial x} \text{ for }f,g\in
   k[[\hbar]][x].
 \end{gather*}
 Note that the $Q$-$\Delta$ lemma does not hold here:
 \begin{gather*}0=\im Q\Delta=\im \Delta Q\neq \ker\Delta\cap \im
   Q=k[x]\oplus k\eta \cap 3x^2k[x]=3x^2 k[x].
 \end{gather*}
 However, we have $Q(f+g\eta)=0\Rightarrow g=0\Rightarrow
 K(f+g\eta)=0.$ Thus, one can take $\beta\colon V\to V[[\hbar]]$
 defined by an inclusion of $Q$ closed elements into $V[[\hbar]]$, and hence, by Theorem
 \ref{main_thm}, $L_\hbar$ is smooth formal.

 \subsection{Physical considerations}
 Differential BV algebras frequently appear in physical situations
 from the quantization of a classical action.  In an odd graded Lie
 algebra, one has the \emph{classical master equation}
 $\left[S^{(0)},S^{(0)}\right]=0$ and a degree zero element $S^{(0)}$
 satisfying $[S^{(0)},S^{(0)}]=0$ is called a \emph{classical action}.
 The operator $Q:V^i \to V^{i+1}$ defined by
 \[Q(v)=[S^{(0)},v]\] satisfies $Q^2=0$ and is a derivation of the Lie
 bracket.  A degree zero element $O^{(0)}\in V$ is called a \emph{classical
   observable} if \[[S^{(0)},O^{(0)}]=0\Leftrightarrow
 Q\left(O^{(0)}\right)=0.\] 
The superscript $(0)$ in $S^{(0)}$ and $O^{(0)}$ is not meant to
indicate degree, it is to emphasize that there are
 no ``$\hbar$'s''  which may appear in quantum actions and quantum observables.

If $(V,\Delta)$ is a BV algebra, then one has the
 \emph{quantum master equation}
 \begin{equation}
   \Delta\left(e^{\frac{S}{\hbar}}\right)=0
 \end{equation}
 and a degree zero solution $S=S^{(0)}+S^{(1)}\hbar
 +S^{(2)}\hbar^2+\cdots \in V[[\hbar]]$ is called a \emph{quantum
   action.}  The quantum master equation
 is equivalent to
 \[\hbar\Delta(S)+\frac{1}{2}[S,S]=0\]
 which unpacks in powers of $\hbar$ as the sequence of equations
 \begin{align*}
   [S^{(0)},S^{(0)}]&=0\\
   \Delta\left(S^{(0)}\right)+[S^{(0)},S^{(1)}]&=0\\
   \Delta\left(S^{(1)}\right)+\frac{1}{2}[S^{(1)},S^{(1)}]&=0\\
   \vdots
 \end{align*}
 So the $S^{(0)}$ term of a quantum action is a classical action.  The
 classical action $S^{(0)}$ is itself a quantum action with no
 additional $\hbar$ terms provided it satisfies
 the additional condition $\Delta\left(S^{(0)}\right)=0$.  The
 condition $\Delta\left(S^{(0)}\right)=0$ also implies that $\Delta
 Q+Q\Delta =0$ and it follows that $(V,Q,\Delta)$ is a differential BV
 algebra. 

 Assume that $S^{(0)}$ is a classical action that is also a
 quantum action so that the triple $(V,Q,\Delta)$ is a differential BV algebra.
 An element $O=O^{(0)}+\hbar O^{(1)}+\hbar^2 O^{(2)}+\cdots
 \in V[[\hbar]]$ called a \emph{quantum observable} if and only if
 \begin{equation*}
   \Delta\left(Oe^{\frac{1}{\hbar}S^{(0)}}\right)=0 \Leftrightarrow K(O)=0.
 \end{equation*}
 The quantum master equation $\Delta \left(
   e^{\frac{S}{\hbar}}\right)=0$ for $S=S^{(0)}+\Gamma$, $\Gamma\in
 V[[\hbar]]$ becomes \[K(\Gamma)+\frac{1}{2}[\Gamma,\Gamma]=0\] which is the
 master equation in the associated quantum dgLa $L_\hbar$.  In
 this physical language, we may summarize our main theorem as:
   \emph{If every classical observable can be extended to a quantum
     observable then there exists a versal quantum action, and
     conversely.}

The significance of having a versal quantum action is that 
the essential physical information contained in quantum correlation
functions can be extracted from this versal action as a
 weak-Frobenius manifold \cite{H} on the space of
classical observables.
See \cite{BK} for the case of Lie
polyvector fields on a Calabi-Yau manifold and \cite{M} for the case of
an abstract dBV algebra satisfying a $Q$-$\Delta$ lemma property.  We
will discuss this physical application and the weak-Frobenius
structure arising
from the smooth formality of $L_\hbar$ in another paper.

\section{Concluding remarks}
Theorem \ref{main_thm} and its proof can be adapted to the more
general context of quantum backgrounds \cite{ptt}.  Quantum
backgrounds provide a kind of generalization of differential BV
algebras where the initial data involves an arbitrary solution to the
quantum master equation---rather than one which has no $\hbar$'s---and
where the operator $\Delta$ may involve components that are
locally polydifferential operators of all orders---rather than being
constrained to at-most second order as in dBV algebras.
In the
language of quantum backgrounds, the result can be stated: \emph{A
  quantum background is smooth formal if and only if every classical
  observable can be extended to a quantum observable, and conversely}.
The very same ideas in the proof of Theorem \ref{main_thm}
for dBV algebras extend to quantum
backgrounds.  Indeed the proof originated in the author's mind in the
quantum background setting.  This paper was written about the
restricted differential BV algebra case just to make it more
accessible and simpler to read.

As a final remark, note that the splitting $\beta$ in the proof of Theorem
\ref{main_thm} gives rise by an automated procedure to a versal
solution to the quantum master equation, the outcome of 
which may be viewed as a
homotopy class of $L_\infty$ maps $H\to L_\hbar$ where $H$ has the
zero $L_\infty$ structure.  By setting $\hbar=0$ one obtains 
a formality map  $H \to L$, but within a special class---not all formality
maps $H\to L$ can originate from a map $\beta:V\to V[[\hbar]]$.
There is a relationship between those
distinguished formality maps that do arise from our splitting $\beta$ 
and the so-called ``special coordinates''
in string theory in which  the quantum correlation functions are equal to
the classical correlation functions.  We will explore this
relationship, as well as the dependence of $\Gamma$ on the
splitting $\beta$ later.

\paragraph{\emph{Acknowledgements}}
  The author is grateful to IH\'ES for providing
  an excellent working environment in January,
  2007.  It was during this period,
  over long discussions and good meals with Jae-Suk
  Park, many of the ideas in this paper were developed.  The author is 
  also grateful to the referee who made helpful comments and excellent
  suggestions which made the paper much better, and to 
  Gabriel Drummond-Cole who read the drafts carefully.

\end{document}